\documentclass[10pt]{amsart}

\usepackage{graphicx}

\newtheorem*{thm}{Theorem}
\newtheorem*{lemma}{Lemma}
\newtheorem{obs}{Observation}
\theoremstyle{definition}

\newcommand{\ignore}[1]{}

\newcommand{\R}{\mathbb{R}}

\newtheorem*{KKM}{The KKM Lemma}
\newtheorem*{LSB}{The LSB Theorem}

\begin{document}
  
\title{The LSB theorem implies the KKM lemma}

\author{Gwen Spencer}
\address{Department of Mathematics\\ Harvey Mudd College\\ Claremont,
  CA  91711}
\email{gspencer@hmc.edu}

\author{Francis Edward Su}
\address{Department of Mathematics\\ Harvey Mudd College\\ Claremont,
  CA  91711}
\email{su@math.hmc.edu}

\thanks{The authors gratefully 
acknowlege partial support by NSF Grant DMS-0301129 (Su).}



\maketitle


Let $S^d$ be the unit $d$-sphere, the set of all
points of unit Euclidean distance from the origin in $\R^{d+1}$. 
Any pair of points in $S^d$ of the form $x,-x$ is a 
pair of {\em antipodes} in $S^d$.
Let $\Delta^d$ be the $d$-simplex formed by the convex hull 
of the standard unit vectors in $\R^{d+1}$.  Equivalently,
$\Delta^d = \{ (x_1,...,x_{d+1}) : \sum_i x_i = 1, x_i\geq 0\}$.
The following are two classical results about closed covers
of these topological spaces:

\begin{LSB}[Lusternik-Schnirelman-Borsuk \cite{LuSc30,Bors33}]
Suppose that $S^d$ is covered by $d+1$ closed sets $A_1,...,A_{d+1}$.
Then some $A_i$ contains a pair of antipodes. 
\end{LSB}

\begin{KKM}[Knaster-Kuratowski-Mazurkiewicz \cite{KKMa29}]
Suppose that $\Delta^d$ is covered by $d+1$ closed sets
$C_1, C_2,...C_{d+1}$ such that 
for each $x$ in $\Delta^d$, $x$ is in $\cup \{ C_i : x_i > 0\}$.  
Then all the sets have a common intersection point, i.e., 
$\cap_{i=1}^{d+1} C_i$ is non-empty.
\end{KKM}

A cover satisfying the condition in the KKM lemma is 
sometimes called a {\em KKM cover}.  It can be rephrased in an
alternate way:
associate labels $1,2,..,d+1$ to the vertices of $\Delta^d$;
then demand that vertex $i$ is covered by set $C_i$
and that each face of $\Delta^d$ is covered by the sets that correspond to the 
vertices spanning that face.  


Both of the above set-covering results are perhaps 
best known in connection with their equivalent formulations in
topology; the LSB theorem is equivalent to the Borsuk-Ulam theorem
\cite{Bors33}, and the KKM lemma is equivalent to the Brouwer fixed
point theorem \cite{KKMa29}.  Also, the LSB theorem has found
spectacular application in proofs of the 
Kneser conjecture in combinatorics \cite{Bara78, Gree02}.
The KKM lemma has numerous applications in economics, e.g., see \cite{Bord85}.

Since the Brouwer fixed point theorem
can be obtained as a consequence of the Borsuk-Ulam theorem
\cite{Su97}, it is natural to ask if the there is a direct proof of the 
KKM lemma using the LSB theorem.
The purpose of this article is to provide such a proof.

\begin{thm}
\label{lsbkkm}
The LSB theorem implies the KKM lemma.
\end{thm}

\ignore{*************************

In order to show that the second follows from the former, we begin by
establishing the following Lemma:

\begin{lemma}
Given any KKM cover of the d-simplex, it is possible to construct a
KKM cover that is ``non-degenerate'' without creating any additional
points in the maximal intersection. Here, ``non-degenerate'' means
that any face of the simplex is covered only by sets corresponding to
the carrying vertices of that face. 
\end{lemma}

\begin{proof}
  Assume some KKM cover of the $d$-simplex by $C_1, C_2, ...C_{d+1}$.
  ``Restrict'' the sets by taking their intersection with the simplex.
  Now the sets Consider the boundary of the simplex crossed with an
  interval with the sets corresponding to the carrier of the face
  extended on the interval...orignal simplex and new part together is
  topologically equivlant to a simplex.  No new maximal intersections
  created. (How to write this?)
\end{proof}

********************************}

\begin{figure}[t]
\begin{center}
\includegraphics[height=2in]{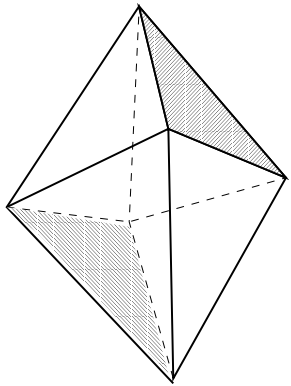}
\end{center}
\caption{The $2$-sphere $\Sigma^2$ in the $L^1$-norm, which is the
  boundary of an octahedron.  The ``top'' and ``bottom'' facets are
  shaded.}
\label{sigmad1}
\end{figure}

Observe that the LSB theorem holds for a $d$-sphere 
under any metric on $\R^{d+1}$, since such a sphere and $S^d$ are
related by an antipode-preserving homeomorphism.  In particular, the
LSB theorem holds for a $d$-sphere under the $L^1$ norm:
$$\Sigma^d := \{ (x_1,...,x_{d+1}) : \sum_i |x_i| = 1 \}.$$
For $d=2$, this is just the boundary of
a regular octahedron, and for general $d$, $\Sigma^d$ is the boundary of 
the {\em $(d+1)$-crosspolytope}.  It is the union of $2^{d+1}$
facets which are simplices, one for each orthant of $\R^{d+1}$.  See Figure
\ref{sigmad1}.



It will be convenient, then, to use the LSB theorem for $\Sigma^d$ 
to prove the KKM lemma, 
because $\Delta^d$ is naturally embedded in $\Sigma^d$; namely,
$\Delta^d$ is the facet of $\Sigma^d$ for which $\sum_i x_i = 1$.  
Call this facet $F_{top}$, the ``top'' facet, and call the antipodal facet
the ``bottom'' facet $F_{bot}$.
Let $F_{mid}$ denote the complement of $F_{top} \cup F_{bot}$ in
$\Sigma^d$, the ``middle'' band of the $d$-sphere.
The strategy of our proof will be to 
assume for the sake of contradiction that a KKM cover of $\Delta^d$
has no common intersection point.  Then we extend these sets to
construct a closed cover of
$\Sigma^d$ whose sets contain no pair of antipodes, thereby 
contradicting the LSB theorem.

\begin{proof}

{\bf Part (I). Construction.}
We first consider the case where a given KKM cover 
$C_1,...,C_{d+1}$ of $\Delta^d$
is {\em non-degenerate}, i.e., for each 
$x$ in $\Delta^d$ and set $C_i$, $x$ is in $C_i$ only if $x_i > 0$. 
In the alternate characterization of the KKM cover,
this means that each face is {\em only} covered by the sets 
that correspond to the vertices spanning that face.  For example, 
the figure at left in Figure \ref{nondeg} is degenerate because the white set
covers a point on the bottom edge of the triangle.

For the sake of contradiction, assume that there is no point common to
all the sets $C_1,...,C_{d+1}$.  
For each $i$, let $-C_i$ be the set in $F_{bot}$ antipodal to $C_i$.
Let $B_i$ be the complement of $-C_i$ in $F_{bot}$.  
By assumption every point of $F_{top}$ is
excluded from at least one $C_i$.  Hence the complementary sets $B_i$
form an open cover of $F_{bot}$ (in the relative topology).
Moreover, the sets $B_i$ satisfy a certain kind of non-degeneracy
that follows from the non-degeneracy of the $C_i$'s:  for $x$ in $F_{bot}$, 
$x_i=0$ implies that $x$ is covered by $B_i$.
By normality, the sets $B_i$ can be shrunk to 
obtain closed subsets $E_i$ of $B_i$ that still cover $F_{bot}$ 
and satisfy the same non-degeneracy.

Now that $F_{bot}$ has been covered, 
we construct a cover of $F_{top}\cup F_{mid}$.
For $x=(x_i)$ in $\Sigma^d$, 
let $\mbox{pos}(x) :=\sum_{x_i>0} x_i$.  Note that $\mbox{pos}(x)=0$ on
$F_{bot}$ but $\mbox{pos}(x)>0$ on $F_{top}$ and $F_{mid}$.
Define a function $f=(f_i)$ on $F_{top} \cup F_{mid}$ by:
$$
f_i(x) =
 \frac{x_i+|x_i|}{2\, \mbox{pos}(x)} \qquad
\textrm{if $\mbox{pos}(x)>0$}
.
$$
Note that $f$ is a continuous function taking $F_{top}\cup F_{mid}$ 
to $F_{top}$, and it fixes $F_{top}$.

Then $D_i := f^{-1}(C_i)$ is a closed subset of $F_{top} \cup F_{mid}$ 
in the relative topology.  
We may think of the set $D_i$ as extending the set 
$C_i$ on $F_{top}$ to cover $F_{mid}$.  In fact, $D_i$ extends
the boundary of $C_i$ in a linear fashion across $F_{mid}$.  
See Figure \ref{lsbsets}.
We record some observations about the sets $D_i$:

\begin{obs}
\label{Dicover}
Since the $C_i$'s cover $F_{top}$, the $D_i$'s cover $F_{top} \cup
F_{mid}$.
\end{obs}

\begin{obs}
\label{Direstrict}
Since $f$ fixes $F_{top}$, each $D_i$ restricted to $F_{top}$ is just $C_i$.
\end{obs}

\begin{obs}
\label{Dixi}
If $x$ is in $D_i$, then $x_i > 0$.
\end{obs}

The first two observations are apparent from the definition of $f$, and
the last observation follows 
by noting that if $x$ is in $D_i$, then $f(x)$ is in $C_i$ and the
non-degeneracy of $C_i$ implies that $f_i(x) > 0$.
But this can only occur if $x_i > 0$.

\begin{figure}[h]
\begin{center}
\includegraphics[height=3in]{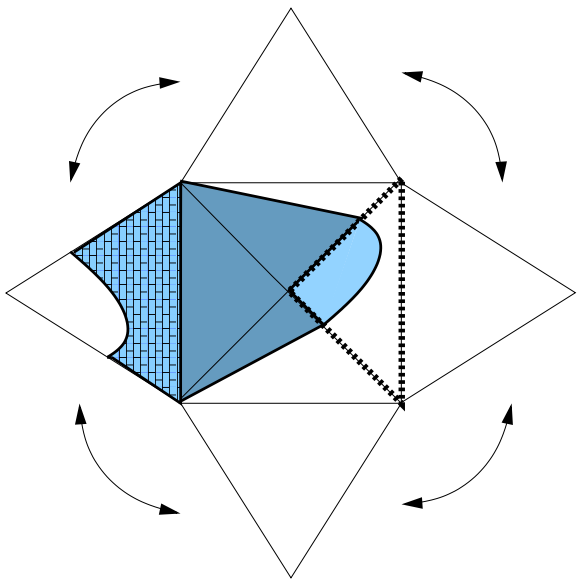}
\end{center}
\caption{The octahedral $2$-sphere $\Sigma^2$ unfolded, with 
shaded set $A_i$ derived from a set $C_i$ in the KKM cover.
The set $A_i$ consists of three regions: light-shaded, dark-shaded,
and bricked.  The light-shaded region is $C_i$; it sits in $F_{top}$
(the triangle with dashed outline).
The set $D_i$ extends $C_i$ and includes both the light-shaded and
dark-shaded regions of $A_i$.  The bricked region is $E_i$; it
sits in the facet antipodal to $F_{top}$.
Note its relation with $C_i$.  
}
\label{lsbsets}
\end{figure}

Now let $A_i=D_i \cup E_i$.  We shall verify that 
the $A_i$'s cover $\Sigma^d$ and are closed sets, yet no $A_i$ 
contains a pair of antipodes.  This verification will 
contradict the LSB theorem, 
forcing us to reject our initial assumption that the
KKM cover had no common intersection point.

{\bf Part (II). Verification.}
Clearly the $A_i$'s cover $\Sigma^d$, by 
Observation \ref{Dicover} and the fact that the $E_i$'s cover $F_{bot}$.

To show $A_i$ is closed, note that $E_i$ is a closed
subset of $\Sigma^d$ and $D_i$ is closed in $F_{top}\cup F_{mid}$ 
(but not necessarily $\Sigma^d$).  Thus 
it suffices to show that any limit points
of $D_i$ in $F_{bot}$ must lie in $E_i$.
Observation \ref{Dixi} implies that 
a limit point $x$ 
of $D_i$ must satisfy $x_i \geq 0$, but since points in $F_{bot}$
have no positive coordinates, 
a limit point of $D_i$ in $F_{bot}$ must satisfy $x_i=0$.
By the non-degeneracy of $E_i$, $x$ must be in $E_i$.

To show that $A_i$ contains no pair of antipodes, we note that $E_i$
cannot contain a pair of antipodes and by Observation \ref{Dixi},
neither can $D_i$, because for any $x$ in $D_i$, 
$x_i$ and $-x_i$ cannot both be positive.  So all that remains is to check
that there is no pair $x$ in $D_i$ and $-x$ in $E_i$.  But this can 
only occur if $x$ is in $F_{top}$.  
By construction $C_i$ cannot have antipodes in $E_i$, so 
Observation \ref{Direstrict} shows that $D_i$ has no antipodes in
$E_i$.  Hence the $A_i$ form a cover of $\Sigma^d$ by $d+1$ closed
sets, yet no $A_i$ contains a pair of antipodes.  This contradicts the
LSB theorem.

\begin{figure}[t]
\begin{center}
\includegraphics[height=2in]{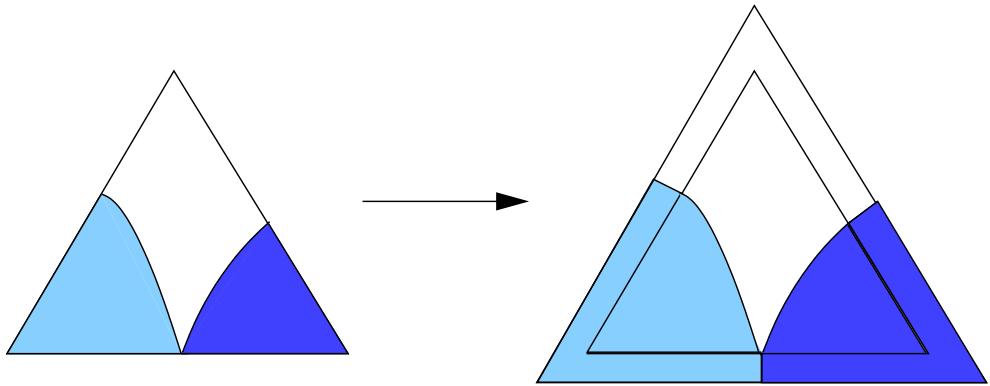}
\end{center}
\caption{In these diagrams, the sets are closed and contain their
 boundaries.  At left, the KKM cover is degenerate because the white
  (non-shaded set) covers a point on the bottom edge.
At right, the same KKM cover has been ``thickened'' to form a
non-degenerate KKM cover.}
\label{nondeg}
\end{figure}

{\bf Part (III). Degenerate KKM covers.}
Finally, we consider the case where the KKM cover is degenerate.  We claim
that a degenerate cover of $\Delta^d$ can be made non-degenerate by
``thickening'' up the boundary and extending the cover in a way that
introduces no new common intersection point.
Let $S$ be the subset of 
$\Delta^d \times [0,1]$ consisting of the points in  
$\Delta^d \times \{0\}$ and $\partial \Delta^d \times [0,1]$.  
(Here $\partial \Delta^d$ denotes the boundary of $\Delta^d$.)
Thus $S$ is homeomorphic to $\Delta^d$; in fact, it is 
$\Delta^d$ with its boundary ``thickened'' up.
Given a KKM cover of $\Delta^d$ by $C_1, C_2,...C_{d+1}$, 
we construct a KKM cover $C'_1,C'_2,...,C'_{d+1}$ of $S$ 
that is non-degenerate.  
First, for $(x,0)$ in $\Delta^d \times \{0\}$, put $(x,0)$ in $C'_i$
if $x$ is in $C_i$.
Then, for $(x,t)$ in $\partial \Delta^d \times [0,1]$ where $t>0$,
put $(x,t)$ in $C'_i$ if $x$ is in $C_i$ {\em and} $x_i>0$.
One may check that the $C'_i$'s are closed and by construction
there are no points of $\cap C'_i$ in the portions of $S$ where $t>0$.
See
Figure \ref{nondeg}.
This ``thickened'', non-degenerate cover can then be used as in the first
part of this proof.
\end{proof}

We remark that although our proof of the KKM lemma
appears non-constructive, the
asserted KKM intersection is hiding in our construction in the
following way.  When we assume (falsely) that the asserted KKM
intersection does not exist, 
we are (wrongly) led to conclude that the $B_i$'s
cover the bottom facet of $\Sigma^d$.  In
actuality, these open sets do not cover the bottom facet;
the set of points that are exposed are precisely the points whose
antipodes comprise the asserted KKM intersection in the top facet.


\end{document}